\documentclass[12TP,draft]{article}

\begin{document}

\begin{center}
\LARGE\noindent\textbf{A Note  on Long non-Hamiltonian Cycles in One Class of Digraphs}\\

\end{center}
\begin{center}
\noindent\textbf{S.Kh. Darbinyan and I.A. Karapetyan}\\

Institute for Informatics and Automation Problems, Armenian National Academy of Sciences

E-mails: samdarbin@ipia.sci.am, isko@ipia.sci.am\\
\end{center}

\textbf{Abstract}\\

 Let $D$ be a strong digraph on $n\geq 4$ vertices. In [3, Discrete Applied Math., 95 (1999) 77-87)], J. Bang-Jensen, Y. Guo and A. Yeo proved the following theorem: if (*) $d(x)+d(y)\geq 2n-1$ and $min \{ d^+(x)+ d^-(y),d^-(x)+ d^+(y)\}\geq n-1$ for every pair of non-adjacent vertices $x, y$ with a common in-neighbour or a common out-neighbour, then $D$ is hamiltonian.
 In this note we show that:  if $D$ is not directed cycle and satisfies the condition (*), then  $D$ contains a cycle of length $n-1$ or $n-2$. \\

Keywords: Digraphs; cycles; Hamiltonian cycles; long non-Hamiltonian cycles\\

\noindent\textbf{1. Introduction and Terminology}\\

We shall assume that the reader is familiar with the standard terminology on directed graphs (digraphs) and refer the reader to monograph of Bang-Jensen and Gutin [1] for terminology not discussed here. In this paper we consider finite digraphs without loops and multiple arcs. For a digraph $D$, we denote by $V(D)$ the vertex set of $D$ and by  $A(D)$ the set of arcs in $D$. Often we will write $D$ instead of $A(D)$ and $V(D)$. The arc of a digraph $D$ directed from $x$ to $y$ is denoted by $xy$. For disjoint subsets $A$ and  $B$ of $V(D)$  we define $A(A\rightarrow B)$ \, as the set $\{xy\in A(D) / x\in A, y\in B\}$ and $A(A,B)=A(A\rightarrow B)\cup A(B\rightarrow A)$. If $x\in V(D)$ and $A=\{x\}$ we write $x$ instead of $\{x\}$.   The out-neighbourhood of a vertex $x$ is the set $N^+(x)=\{y\in V(D) / xy\in A(D)\}$ and $N^-(x)=\{y\in V(D) / yx\in A(D)\}$ is the in-neighbourhood of $x$. Similarly, if $A\subseteq V(D)$ then $N^+(x,A)=\{y\in A / xy\in A(D)\}$ and $N^-(x,A)=\{y\in A / yx\in A(D)\}$. We call the vertices in $N^+(x)$, $N^-(x)$, the out-neighbours and in-neighbours of $x$. The out-degree of $x$ is $d^+(x)=|N^+(x)|$ and $d^-(x)=|N^-(x)|$ is the in-degree of $x$. The out-degree and in-degree of $x$ we call its semi-degrees. Similarly, $d^+(x,A)=|N^+(x,A)|$ and $d^-(x,A)=|N^-(x,A)|$. The degree of the vertex $x$ in $D$ defined as $d(x)=d^+(x)+d^-(x)$ (similarly, $d(x,A)=d^+(A)+d^-(A)$). The subdigraph of $D$ induced by a subset $A$ of $V(D)$ is denoted by $\langle A\rangle$. The path (respectively, the cycle) consisting of the distinct vertices $x_1,x_2,\ldots ,x_m$ ( $m\geq 2 $) and the arcs $x_ix_{i+1}$, $i\in [1,m-1]$  (respectively, $x_ix_{i+1}$, $i\in [1,m-1]$, and $x_mx_1$), is denoted  $x_1x_2\cdots x_m$ (respectively, $x_1x_2\cdots x_mx_1$). For a cycle  $C_k=x_1x_2\cdots x_kx_1$, the subscripts considered modulo $k$, i.e. $x_s=x_i$ for every $s$ and $i$ such that  $i\equiv s\, (\hbox {mod} \,k)$. If $P$ is a path containing a subpath from $x$ to $y$ we let $P[x,y]$ denote that subpath. Similarly, if $C$ is a cycle containing vertices $x$ and $y$, $C[x,y]$ denotes the subpath of $C$ from $x$ to $y$.
 A digraph $D$ is strongly connected (or just strong) if there exists a path from $x$ to $y$ and a path from $y$ to $x$ in $D$ for every choice of distinct vertices $x$,\,$y$ of $D$. We will denote the  complete bipartite digraph  with partite sets of cardinalities $p$, $q$ by $K^*_{p,q}$. Two distinct vertices $x$ and $y$ are adjacent if $xy\in A(D)$ or $yx\in A(D) $ (or both). We denote by $a(x,y)$ the number of arcs between the vertices $x$ and $y$. In particular, $a(x,y)=0$ (respectively, $a(x,y)\not=0$) means that $x$ and $y$ are not adjacent (respectively, are adjacent). 

For integers $a$ and $b$, $a\leq b$, let $[a,b]$  denote the set of all integers which are not less than $a$ and are not greater than $b$. The digraph $D$ is hamiltonian (is pancyclic, respectively) if it contains a hamiltonian cycle, i.e. a cycle of length $|V(D)|$ (contains a cycle of length $m$ for any $3\leq m \leq |V(D)|$).\\

Meyniel [12] proved the following theorem: if $D$ is a strong digraph on $n\geq 2$ vertices and $d(x)+d(y)\geq 2n-1$ for all pairs of non-adjacent vertices in $D$, then $D$ is hamiltonian (for short proofs of Meyniel's theorem see [4, 13]). 

Thomassen [15] (for $n=2k+1$) and Darbinyan [6] (for $n=2k$) proved: if $D$ is a digraph on $n\geq 5$ vertices with minimum degree at least $n-1$ and with minimum semi-degree at least $n/2-1$, then $D$ is hamiltonian (unless some extremal cases). \\

In each above mentioned theorems (as well as, in theorems Ghouila-Houri [10], Woodall [16], Manoussakis [11]) imposes a degree condition on all pairs of non-adjacent vertices (on all vertices). Bang-Jensen, Gutin, Li, Guo and Yeo [2, 3] obtained sufficient conditions for hamiltonicity of digraphs in which degree conditions requiring only for some pairs of non-adjacent vertices. Namely, they proved the following theorems (in all three theorems  $D$ is a strong digraph on $n\geq 2$ vertices).\\

\noindent\textbf{Theorem A} [1, 2]. If $min \{d(x),d(y)\}\geq n-1$ and $d(x)+d(y)\geq 2n-1$ for every pair of non-adjacent vertices  $x$, $y$ with a common in-neighbour, then $D$ is hamiltonian.

\noindent\textbf{Theorem B} [1, 2]. If $min \{d^+(x)+d^-(y),d^-(x)+d^+(y)\}\geq n$ for every pair of non-adjacent vertices $x$, $y$ with a common out-neighbour or a common in-neighbour, then $D$ is hamiltonian.

\noindent\textbf{Theorem C} [3]. If $min \{d^+(x)+d^-(y),d^-(x)+d^+(y)\}\geq n-1$ and $d(x)+d(y)\geq 2n-1$ for every pair of non-adjacent vertices $x$, $y$ with a common out-neighbour or a common in-neighbour, then $D$ is hamiltonian. \\

 Note that  Theorem C generalizes Theorem B.
In [9, 14, 5, 7] it was shown  that if the strong digraph $D$  satisfies the condition of the theorem of Ghouila-Houri [10] (Woodall [16], Meyniel [12], Thomassen and Darbinyan [15, 6]), then $D$ is pancyclic (unless some extremal cases, which are characterized). In [8], we possed  the following problem:\\

\noindent\textbf{Problem}. Characterize those digraphs which satisfy the conditions of Theorem A (B, C), but are not pancyclic.\\

In [8], we have show that: 

(i) if a strong digraph $D$ satisfies the condition of Theorem A and the minimum semi-degree of   $D$ at least two; or 

(ii) if a strong digraph  $D$ is not directed cycle and satisfies the condition  of Theorem B, then either $D$ contains a cycle of length $n-1$ or $n$ is even and $D$ is isomorphic to complete bipartite digraph or to complete bipartite digraph minus one arc.\\ 

In [8], we also possed  the following 

\noindent\textbf {Conjecture}. Let a digraph $D$ on $n\geq 4$ vertices satisfies the conditions of Theorem C. Then $D$ contains a cycle of length $n-1$ maybe except some digraphs which has a "simple" characterization.\\

Support for the our conjecture, in this note by using the proof of Theorem C ( Theorem 3.1, [3]), we showe that: if $D$ is not directed cycle and satisfies the conditions of Theorem C, then  $D$ contains a cycle of length $n-1$ or $n-2$. \\

\noindent\textbf{2. Preliminaries }\\

The following well-known simple lemmas is the basis of our results and other theorems on directed cycles and paths in digraphs. It we  will be used extensively in the proof of our result.\\

\noindent\textbf{Lemma 1} [9]. Let $D$ be a digraph on $n\geq 3$ vertices containing a cycle $C_m$, $m\in [2,n-1] $. Let $x$ be a vertex not contained in this cycle. If $d(x,C_m)\geq m+1$, then  $D$ contains a cycle $C_k$ for all  $k\in [2,m+1]$.   \\

\noindent\textbf{Lemma 2} [4]. Let $D$ be a digraph on $n\geq 3$ vertices containing a path $P:=x_1x_2\ldots x_m$, $m\in [2,n-1]$ and let $x$ be a vertex not contained in this path. If one of the following conditions holds:

 (i) $d(x,P)\geq m+2$; 

 (ii) $d(x,P)\geq m+1$ and $xx_1\notin D$ or $x_mx_1\notin D$; 

 (iii) $d(x,P)\geq m$, $xx_1\notin D$ and $x_mx\notin D$;

\noindent\textbf{}then there is an  $i\in [1,m-1]$ such that $x_ix,xx_{i+1}\in D$, i.e. $D$ contains a path $x_1x_2\ldots x_ixx_{i+1}\ldots x_m$ of length $m$  (we say that  $x$ can be inserted into $P$ or the arc  $x_ix_{i+1}$ is a partner of $x$ on $P$). \\

\noindent\textbf{Lemma 3 }[2]. Let $P:=x_1x_2\ldots x_m$  be a path in $D$ and let $x$, $y$ be vertices of $V(D)-V(P)$ (possibly $x=y$). If there do not exist consecutive vertices $x_i, x_{i+1}$ on $P$ such that $x_ix$, $yx_{i+1}$ are arcs of $D$, then $d^-(x,P)+d^+(y,P)\leq m+\epsilon$, where $\epsilon =1$ if $x_mx\in D$ and 0, otherwise. \\

\noindent\textbf{3. Main result}\\

Let $C$ be a cycle in digraph $D$. For the cycle $C$, a $C$-bypass is an $(x,y)$-path $P$ of length at least two with both end-vertices $x$ and $y$ on $C$ and no other vertices on $C$. The length of the path $C[x,y]$ is the gap of $P$ with respect to $C$.

 If $\{x,y\}$ is a pair of non-adjacent vertices with a common in-neighbour or a common out-neighbour, then in the proof of the theorem we say that $\{x,y\}$ is a good pair.\\

In the proof of the our theorem we use (in the main) the notations which are used in the proof of Theorem C (Theorem 3.1, [3]).\\

\noindent\textbf{Theorem}. Let $D$ be a strong digraph with $n\geq 2$ vertices, which is not directed cycle. Suppose that $min \{d^+(x)+d^-(y),d^-(x)+d^+(y)\}\geq n-1$ and $d(x)+d(y)\geq 2n-1$ for every pair of non-adjacent vertices $x$, $y$ with a common out-neighbour or a common in-neighbour, then $D$ contains a cycle of length $n-2$ or $n-1$. \\

\noindent\textbf{Proof}. Suppose, to the contrary, that $D$ contain no cycles of length $n-2$ or $n-1$. Let $C:=x_1x_2\ldots x_mx_1$ be a longest non-hamiltonian cycle in $D$. Then  $3\leq m\leq n-3$ and let $R:=V(D)-V(C)$. Observe that if $y\notin V(C)$, then $y$ has no partner on $C$. We shall use this often without explicit reference. For the digraph $D$ provided that $D$ is not hamiltonian, in [3] (Theorem 3.1), J.Bang-Jensen, Y. Guo and  A. Yeo proved the following Claims 1 and 2.\\

\noindent\textbf{Claim 1}. Let $y$ be a vertex of $R$. If $x_{\alpha}\not =x_{\beta}$, $x_{\alpha}y$, $yx_{\beta} \in D$ and $A(y,C \setminus V(C[x_{\beta},x_{\alpha}])=\emptyset$, then the following holds:
$$|V(C^{'})|\geq 1, \quad d(y,C)=d(y,C^{''})=|C^{''}|+1,  \eqno (1) $$
$$ d^+(x_{\beta -1},C^{''})+d^-(x_{\alpha +1},C^{''})=|C^{''}|+1,  \eqno (2) $$
$$ d(y,R)+ d^+(x_{\beta -1},R)+d^-(x_{\alpha +1},R)=2(n-m-1),  \eqno (3) $$
$$  d^+(x_{\beta -1},C^{'})=d^-(x_{\alpha +1},C^{'})=|C^{'}|-1,  \eqno (4) $$
where $C^{'}:=C[x_{\alpha +1},x_{\beta -1}]$ and $C^{''}:=C[x_{\beta},x_{\alpha}]$.    \framebox  \\\\ 

\noindent\textbf{Claim 2}. $D$ contains a $C$-bypass.  \framebox  \\\\

Note that Claims 1 and 2 also are true if in $D$ a longest non-hamiltonian cycle has length at most $n-3$ (the proofs are just the same). 

From (4) it follows that if $|C^{'}|\geq 2$, then $P:=x_{\beta -1}x_{\alpha +2} \ldots x _{\beta -2}x_{\alpha +1}$ is a hamiltonian $(x_{\beta -1},x_{\alpha +1})$-path in $\langle C^{'}\rangle$. Therefore, similarly (2), we obtain (Lemma 3)
$$ d^-(x_{\beta -1},C^{''})+d^+(x_{\alpha +1},C^{''})=|C^{''}|+1.  $$
Combining this last inequality with (2) yields
$$ d(x_{\beta -1},C^{''})+d(x_{\alpha +1},C^{''})\leq 2|C^{''}|+2.  \eqno (5) $$

We now  prove the following claim: 

\noindent\textbf{Claim 3}. Let $x_myx_{\gamma +1}$ be a $C$-bypass and $A(y,C[x_1,x_{\gamma}])=\emptyset$. Then $\gamma \geq 3$.

\noindent\textbf{Proof}. Suppose that $ \gamma \leq 2$. Let now $C^{''}:=C[x_{\gamma +1},x_m]$. We shall consider the cases $\gamma =1$, $\gamma =2$ separately.

\noindent\textbf{Case 1}. $\gamma =1$. Then similarly (1) and (3) we have $d(y,C)=d(y,C^{''})\leq m$, $d(x_1,C^{''})\leq m $ and $d(y,R)+d(x_1,R)\leq 2(n-m-1)$. Therefore, since $\{y,x_1\}$ is a good pair and $|C^{''}|=m-1$, we have 
$$
2n-1\leq d(y)+d(x_1)=d(y,R\cup C^{''})+d(x_1,R\cup C^{''})\leq 2(n-m-1)+2|C^{''}|+2=2n-2,
$$ a contradiction.

\noindent\textbf{Case 2}. $\gamma =2$. Then, since $|R|\geq 3$, for any $i\in [1,2]$ we obtain that $d(y,R)+d(x_i,R)\leq 2(n-m-1)$. Since $\{y,x_i\}$ is a good pair and (1), it follows that 
$$
2n-1\leq d(y)+d(x_i)\leq 2(n-m-1)+d(y,C^{''})+d(x_i,C^{''})+2 \leq 2(n-m-1)+|C^{''}|+3+d(x_i,C^{''}).
$$ 
From this we obtain that $d(x_i,C^{''})\geq m=|C^{''}|+2$. Hence, by Lemma 2, the vertex $x_1$  ($x_2$) has a partner on $C^{''}$. Therefore there is a $(x_{3},x_m)$-path with vertex  set $V(C)$. This path with the vertex $y$ forms a  non-hamiltonian cycle longer than $C$. Claim 3 is proved.  \framebox  \\\\

Let  $P:=u_1u_2\ldots u_s$ ($s\geq 3$) be a $C$-bypass with minimum gap among the gaps of all $C$-bypasses. Assume w.l.o.g. that $P$ is minimal with respect to the minimum gap and let $u_1:=x_1$, $u_s:=x_{\gamma }$ with $2\leq \gamma \leq m$. 

 In the following we suppose, further, that $\gamma =2$ (the proof for the case $\gamma \geq 3$ is same as the proof of Theorem C (Theorem 3.1, [3])).

 Then $R=\{u_2,u_3, \ldots , u_{s-1}\}$, $s\geq 5$ and for any pair of  $i,j$ with $2\leq i< j \leq s-1$
$$
u_iu_j\in D \quad \hbox{if and only if}\quad j=i+1.   \eqno (6)
$$
Since $|R|\geq 3$ and $C$ is a longest non-hamiltonian cycle in $D$, it is easy to see that 
$$
x_1u_{s-1}\notin D, \, u_2x_2\notin D \quad \hbox {and} \quad 
d^-(u_2,\{x_{m-1},x_m\})=d^+(u_{s-1},\{x_{3},x_4\})=0.    \eqno (7)
$$

\noindent\textbf{Case 1}. $x_2u_2\notin D$ and there is an $i\in [3,m]$ such that $x_iu_2\in D$. Then by (7) we have $a(x_2,u_2)=0$, $3\leq i\leq m-2$ and by Claim 3, $d^+(u_2,\{x_3,x_4\})=0$. 

Assume that $d^+(u_2,C[x_{i+1},x_1])\not=0$. Then there are integers $l$ and $j$ with $i\leq l\leq j-1\leq m$ such that $x_lu_2, u_2x_j\in D$ and $A(u_2,C[x_{l+1},x_{j-1}])=\emptyset$. By (1), 
$$d(u_2,C)=d(u_2,C[x_{j},x_{l}])= |C[x_{j},x_{l}]|+1.$$
On the other hand, since $u_2x_3\notin D$, using Lemma 2 we obtain that
$$
d(u_2,C)=d(u_2,C[x_{3},x_{l}])+d(u_2,C[x_j,x_1])\leq |C[x_{3},x_{l}]|+|C[x_{j},x_{1}]|+1=|C[x_{j},x_{l}]|,
$$
a contradiction.

Now assume that $d^+(u_2,C[x_{i+1},x_{1}])=0$. Let $i$ is minimal as possible, i.e. $d^-(u_2,C[x_{2},x_{i-1}])=0$. Then by (7) we have $A(u_2,\{x_{m-1},x_m\})=\emptyset$. Let $x_ju_2\in D$, $i\leq j\leq m-2$ and let $j$ is maximal with these properties. If $d^+(u_2,C[x_{2},x_{i}])=0$, i.e. $d^+(u_2,C)=0$, then $d^+(u_2)=1$ because of (6). Since $\{u_2,x_{j+1}\}$ is a good pair, by the condition of the theorem we have $d^+(u_2)+d^-(x_{j+1})\geq n-1$. Therefore $d^-(x_{j+1})\geq n-2$. On the other hand, it is easy to check that  $d^-(x_{j+1},\{u_2,u_3\})=0$, and hence,  $d^-(x_{j+1})\leq n-3$, a contradiction. So we can assume that $d^+(u_2,C[x_{2},x_{i}])\not=0$. Let $u_2x_k\in D$, where $x_k\in C[x_2,x_i]$, and $k$ is minimal as possible. Then, from the minimality of $i$ and $k$ it follows that $A(u_2,C[x_2,x_{k-1}])=\emptyset$. Hence, by Claim 3, $k\geq 5$. By (1) (Claim 1) we have
$$
d(u_2,C)=d(u_2,C[x_k,x_{1}])=|C[x_k,x_{1}])|+1.
$$
On the other hand, since $A(u_2,\{x_{m-1},x_{m}\})=\emptyset$ and $u_2x_1\notin D$, using Lemma 2 we obtain that 
$$
d(u_2,C)=d(u_2,C[x_k,x_{m-2}])+a(u_2,x_1)\leq |C[x_k,x_{m-2}]|+2\leq |C[x_k,x_{1}]|,
$$
a contradiction.\\

\noindent\textbf{Case 2}. $x_2u_2\notin D$ and  $x_iu_2\notin D$ for every $i\in [3,m]$. Then $d^-(u_2,C[x_{2},x_{m}])=0$. 

First assume that there is a $x_i\in C[x_2,x_m]$ such that $u_2x_i\in D$ and $A(u_2,C[x_2,x_{i-1}])=\emptyset$. By Claim 3, $i\geq 5$. Let now $C^{''}:=C[x_i,x_1]$ and $C^{'}:=C[x_2,x_{i-1}]$. Note that $|C^{''}|+|C^{'}|=m$.\\

 \noindent\textbf{Remark}. It is a simple matter to check  that (i) $d(x_{i-1},C{'})\leq 2|C^{'}|-3$  since $x_{i-3}x_{i-1}\notin D$;  (ii) if $|C^{'}|=3$, then $d(x_{2},C{'})\leq 2|C^{'}|-3$ since $x_2x_4\notin D$ and (iii) if $C^{'}|\geq 4$, then $d(x_{2},C{'})\leq 2|C^{'}|-4$ since $d^+(x_2,\{x_4,x_5\})=0$.

By (5) we have  that 
$$ 
d(x_{2},C^{''})+d(x_{i-1},C^{''})\leq 2|C^{''}|+2.  \eqno (8) 
$$
It is not difficult to see that
$$
d^-(x_{2},R)=1, \quad d^+(x_{2},R)\leq n-m-1 \quad \hbox{and hence,} \quad d(x_{2},R)\leq n-m.  \eqno (9) 
$$
$$
d^-(x_{i-1},R)\leq n-m-1, \,\, d^+(x_{i-1},R)\leq n-m-1-l \,\,,\hbox{where} \,\, l:=d^-(u_2,R), \,\, \hbox{and} \,\, d(u_2,R)\leq l+1. \eqno (10)
$$
From (8) it follows that $d(x_{2},C^{''})\leq |C^{''}|+1$ or $d(x_{i-1},C^{''})\leq |C^{''}|+1$. 

Let $d(x_{2},C^{''})\leq |C^{''}|+1$. Then, since $\{x_2,u_2\}$ is a good pair, by (1), (6) and (9) we have
$$
2n-1\leq d(u_2)+d(x_2)=d(u_2,R\cup C^{''})+d(x_2,R\cup C^{''}\cup C{'})\leq 2n-2m+2+2|C^{''}|+d(x_2,C{'}).
$$
Hence $d(x_2,C{'})\geq 2|C^{'}|-3$,  $d(u_2,R)= n-m$ and $d^+(x_2,R)=n-m-1$. By Remark, $|C^{'}|=3$ (i.e., $i=5$). Therefore $u_{s-1}u_2,x_2u_3\in D$ and hence, $x_2u_3\ldots u_{s-1}u_2x_5\ldots x_mx_1x_2$ is a cycle of length $n-2$, which is a contradiction.

Let now $d(x_{i-1},C^{''})\leq |C^{''}|+1$. Since $\{u_2,x_{i-1}\}$ is a good pair, $|C^{'}|+|C^{''}|=m$ and $x_{i-3}x_{i-1}\notin D$, using (1), (10) and Remark we obtain,
$$
2n-1\leq d(u_2)+d(x_{i-1})=d(u_2,R\cup C^{''})+d^+(x_{i-1},R)+d^-(x_{i-1},R)+d(x_{i-1},C^{''}\cup C{'})\leq 2n-2,
$$
a contradiction. 

Second assume that $d^+(u_2,C[x_{2},x_{m}])=0$. Then $A(u_2,C[x_2,x_{m}])=\emptyset$. Since $\{x_2,u_2\}$ is a good pair and $d(x_2,R)\leq n-m$, this implies that 
$$
2n-1\leq d(u_2)+d(x_{2})\leq 2+d(u_2,R)+d(x_2,C\cup R)\leq 2+2n-2m+d(x_{2},C).
$$
From this we obtain that $d(x_2,C)=2m-3$ since $x_2x_4\notin D$ (recall that $m\geq 3$). Now it is not difficult to see that $m=3$, $x_2x_1\notin D$, $u_2x_1,u_{s-1}u_2, x_2u_3\in D$ and  $x_2u_3\ldots u_{s-1}u_2x_1x_2$ is a cycle of length $n-1$, a contradiction.

\noindent\textbf{Case 3}. $x_2u_2\in D$. We can assume that $u_{s-1}x_1\in D$ (otherwise in the converse digraph of D we will have the considered Case 1 or 2).

First assume that $d^+(u_2,C[x_{3},x_{1}])\not=0$. Choose $x_i$ so that $u_2x_i\in D$ and $|C[x_{2},x_{i-1}]|$ is as small as possible. Since $x_2u_2\in D$, there is a $x_j\in C[x_{2},x_{i-1}]$ such that $x_ju_2\in D$ and $A(u_2,C[x_{j+1},x_{i-1}])=\emptyset$. Let now $C^{''}:=C[x_{i},x_j]$ and  $C^{'}:=C[x_{j+1},x_{i-1}]$. By Claim 3 and (1) (Claim 1) we have 
$$
|C^{'}|\geq 3 \quad \hbox{and} \quad d(u_2,C)=d(u_2,C^{''})=|C^{''}|+1. \eqno (11)
$$
Since $x_mu_2\notin D$, using Claim 3, (11) and Lemma 2 it is not difficult to obtain that 
$$
N^-(u_2,C)=\{x_1,x_2,\ldots, x_j\} \quad \hbox{and} \quad N^+(u_2,C)= \{x_i,x_{i+1},\ldots,x_m, x_1\}. \eqno (12)
$$
If $u_{s-1}u_2\notin D$, then $u_{s-1}$ and $u_2$ are not adjacent and hence, $\{u_2,u_{s-1}\}$ is a good pair since $u_2x_1, u_{s-1}x_1 \in D$. Now from (6), (11) and the condition of the theorem it follows that 
$$d(u_2)=d(u_2,C)+d(u_2,R)\leq n-m-1+m-2=n-3$$
and 
$$ n+2\leq d(u_{s-1})=d(u_{s-1},R)+d(u_{s-1},C)\leq n-m-1 + d(u_{s-1},C).$$
Therefore $d(u_{s-1},C)\geq m+3$, and by Lemma 1 $u_{s-1} $ has a partner on $C$, which is a contradiction. So we can assume that this is not the case, i.e. $u_{s-1}u_2\in D$. Then by (12) and the maximality of the cycle $C$ we conclude that 
$$ 
d^+(x_{j+1},R)\leq n-m-1, \,\,d^-(x_{j+1},R)=0,  \,\,\hbox {and hence,}\,\, d(x_{j+1},R)\leq n-m-1. \eqno (14)$$
$$ 
d^-(x_{i-1},R)\leq n-m-1, \,\,d^+(x_{i-1},R)=d^+(x_{i-1},\{u_3\}). \eqno (15)$$
If $d^+(x_{i-1},\{u_3\})=1$, then $x_{i-1}u_3\in D$ and $u_2x_{i+1}\notin D$. Therefore by (12), $x_i=x_1$ and $x_{i-1}=x_m$. Hence we have $x_mu_3\in D$ and the cycle $x_mu_3\ldots u_{s-1}x_2\ldots x_m$ longer than $C$, which is a contradition. So we can assume that $d^+(x_{i-1},R)=0$ and therefore, $d(x_{i-1},R)\leq n-m-1$.  

From (15) it follows that $d(x_{i-1},C^{''})\leq |C^{''}|+1$ or $d(x_{j+1},C^{''})\leq |C^{''}|+1$. Assume that $d(x_{i-1},C^{''})\leq |C^{''}|+1$. Then, since $\{u_2,x_{i-1}\}$ is a good pair, by (1), (14) and $x_{i-3}x_{i-1}\notin D$ we have 
$$
2n-1\leq d(u_2)+d(x_{i-1})=d(u_2,R\cup C)+d(x_{i-1},R)+d(x_{i-1},C^{''}\cup C{'})\leq 2n-2,
$$
a contradiction. Similarly we obtain a contradiction if we assume that $d(x_{j+1},C^{''})\leq |C^{''}|+1$. 

Second assume that $d^+(u_2,C[x_{3},x_{1}])=0$. From $x_2u_2\in D$ and $d^-(u_2,\{x_{m-1},x_m\})=0$ it follows that there is a $x_j\in C[x_2,x_{m-2}]$ such that  $x_ju_2\in D$ and   $A(u_2,C[x_{j+1},x_m ])=\emptyset$. Note that  $\{u_2,x_{j+1}\}$ is a good pair. Then  $d^-(x_{j+1})\geq n-2$ since $d^+(u_2)=1$. From $d^-(x_{j+1},R)\leq 1$ implies that $d^-(x_{j+1},C)\geq n-3\geq m$, which is impossible. Hence in all possible cases we reach a contradiction. The proof of the theorem is complete.

\noindent\textbf {References}\\

[1] J. Bang-Jensen, G. Gutin, Digraphs: Theory, Algorithms and Applications, Springer, 2000.

[2] J. Bang-Jensen, G. Gutin, H. Li, Sufficient conditions for a digraph to be hamiltonian, J. Graph Theory 22 (2) (1996) 181-187.

[3] J. Bang-Jensen, Y. Guo, A.Yeo, A new sufficient condition for a digraph to be hamiltonian, Discrete Applied Math., 95 (1999) 77-87. 

[4] J.A. Bondy, C. Thomassen, A short proof of Meyniel's theorem, Discrete Math. 19 (1977) 195-197.

[5] S.Kh. Darbinyan,  Pancyclic and panconnected digraphs, Ph. D. Thesis, Institute  Mathematici Akad. Nauk BSSR, Minsk, 1981 (see also, Pancyclicity of digraphs with the Meyniel condition, Studia Sci. Math. Hungar., 20 (1-4) (1985) 95-117, in Russian).

[6] S.Kh. Darbinyan, A sufficient condition for the Hamiltonian property of digraphs with  large semidegrees, Akad. Nauk Armyan. SSR Dokl. 82 (1) (1986) 6-8 (see also, arXiv: 1111.1843v1 [math.CO] 8 Nov 2011).

[7] S.Kh. Darbinyan, On the  pancyclicity of digraphs with large semidegrees,  Akad. Nauk Armyan. SSR Dokl. 83 (3) (1986) 99-101 (see also, arXiv: 1111.1841v1 [math.CO] 8 Nov 2011).

[8] S.Kh. Darbinyan, I.A. Karapetyan, On longest non-hamiltonian cycles in digraphs with the conditions of Bang-Jensen, Gutin and Li, Preprint available at http:// arXiv.org/abs/1207.5643v1 [math. CO] 24 Jul 2012. 

[9] R. H\"{a}ggkvist, C. Thomassen, On pancyclic digraphs, J. Combin. Theory Ser. B 20 (1976) 20-40.

[10] A. Ghouila-Houri, Une condition suffisante d'existence d'un circuit hamiltonien, C. R. Acad. Sci. Paris Ser. A-B 251 (1960) 495-497.

[11]  Y. Manoussakis, Directed Hamiltonian graphs, J. Graph Theory 16 (1992) 51-59.

[12] M. Meyniel, Une condition suffisante d'existence d'un circuit hamiltonien dans un graphe oriente, J. Combin. Theory Ser. B 14 (1973) 137-147.

[13] M. Overbeck-Larisch, Hamiltonian paths in oriented graphs, J. Combin. Theory Ser. B 21 (1) (1976) 76-80.

[14] C. Thomassen, An Ore-type condition implying a digraph to be pancyclic, Discrete Math. 19 (1) (1977) 85-92.

[15] C. Thomassen, Long cycles in digraphs,  Proc. London Math. Soc. (3) 42 (1981) 231-251.

[16] D.R. Woodall, Sufficient conditions for circuits in graphs, Proc. London Math. Soc. 24 (1972) 739-755.

\end{document}